\documentclass[letterpaper,11pt]{article}

\usepackage{lineno,hyperref}
\usepackage{amssymb}
\usepackage{color}
\usepackage{cancel}
\usepackage[left=2.8cm,top=2.8cm,right=3.4cm,bottom=4cm]{geometry}
\usepackage{stmaryrd}
\usepackage{amsmath,mathrsfs,amssymb,amsfonts,nicefrac}

\hypersetup{linkcolor=red,citecolor=blue}

\modulolinenumbers[250]
\newcommand{\RR}{ I\!\!R}
\newcommand{\C}{\mathbb{C}}
\newcommand{\R}{\mathbb{R}}

\newcommand{\N}{\mathbb{N}}

\newcommand{\beq}{\begin{equation} }
\newcommand{\eqq}{\end{equation} }
\newcommand{\cuad}{{\sqcap\kern-.68em\sqcup}}
\newcommand{\norm}[1]{\left\|#1\right\|}

\renewcommand{\a}{\alpha}
\newcommand{\eps}{\varepsilon}
\newcommand{\abs}[1]{\left|#1\right|}

\newtheorem{teo}{Theorem}[section]

\newtheorem{lema}[teo]{Lemma}

\newcommand{\bremark}{\begin{remark} \em}
\newcommand{\eremark}{\end{remark} }

\newcommand{\re}{\Re\mathrm{e}}

\def\beeq{\begin{equation}}
\def\eeq{\end{equation}}
\newcommand{\begeqaet}{\begin{eqnarray*}}
\newcommand{\eneqaet}{\end{eqnarray*}}

\begin{document}
\ \\ \\ \\
\begin{center}{\bf\Large Exponential propagation for fractional reaction-diffusion cooperative systems with fast decaying initial conditions.}\medskip

\bigskip

\bigskip

{Anne-Charline COULON$^{a}$ and Miguel YANGARI$^{a,b}$}\\

\ \\$^{a}$\small{Institut de Math\'ematiques, Universit\'e Paul Sabatier \\
118 Route de Narbonne,  F-31062 Toulouse Cedex 4, France.}\\

\medskip
$^{b}$\small{Departamento de Ingenier\'{\i}a  Matem\'atica, Universidad de Chile\\
Blanco Encalada 2120, Santiago, Chile.}

\end{center}

\medskip

\begin{abstract}
We study the time asymptotic propagation of sectorial solutions to the fractional reaction-diffusion cooperative systems. We prove
that the propagation speed is exponential in time, and we find the precise exponent of propagation. This exponent depends on the smallest index of the fractional laplacians and on the principal eigenvalue of the matrix $DF(0)$ where $F$ is the reaction term. We also note
that this speed does not depend on the space direction.
\end{abstract}
\footnotetext[1]{\emph{Email addresses:} anne-charline.coulon@math.univ-toulouse.fr, myangari@dim.uchile.cl}
\footnotetext[2]{\emph{Acknowledgements:} The research leading to these results has received funding from the European Research Council under the European Unions Seventh Framework Programme (FP/2007-2013) / ERC Grant Agreement n.321186 - ReaDi -Reaction-Diffusion Equations, Propagation and Modelling. M. Y. was supported by Becas de Doctorado, Conicyt-Chile and Senescyt-Ecuador. The authors thank Professor J.-M. Roquejoffre for fruitful discussions.}

\date{}
\setcounter{equation}{0}
\section{Introduction}
The reaction diffusion equation with Fisher-KPP nonlinearity \begin{equation}\label{1}\partial_t u +(-\triangle)^{\alpha}u=f(u)\end{equation}
with $\alpha=1$, has been the subject of intense research since the seminal work by Kolmogorov, Petrovskii, and Piskunov \cite{17}.
Of particular interest are the results  of Aronson and Weinberger \cite{1} which describe the evolution of compactly supported data.
They showed that there exists a critical threshold $c^*=2\sqrt{f'(0)}$ such that, for any compactly supported initial value $u_{0}$ in $[0,1]$, if $c>c^*$ then  $u(t, x)\to  0$ uniformly in $\{ |x| \ge ct\}$  as $t\to+\infty$ and if $c<c^*$ then
$u(t, x)\to  1$ uniformly in $\{ |x| \le ct\}$  as $t\to+\infty$. This corresponds to a linear propagation of the fronts.
In addition, (\ref{1}) admits planar travelling wave solutions connecting $0$ and $1$.

Reaction-diffusion equations with fractional Laplacian, that is when $\alpha\in(0,1)$ in \eqref{1}, appear in physical models when the diffusive phenomena are better described by L\'evy processes allowing long jumps, than by Brownian processes - obtained when $\a=1$. The L\'evy processes occur widely in physics, chemistry and biology. Recently these models have attracted much interest. In connection with the discussion given above, in the recent paper \cite{6},  Cabr\'e and Roquejoffre showed that for any compactly supported initial condition, or more generally for initial values decaying faster than $|x|^{-d-2\alpha}$, where $d$ is the dimension of the spatial variable, the speed of propagation becomes exponential in time. They also showed that no travelling wave exist. Their result was sharpened and extended in \cite{7}, who proposed a new (and more flexible) argument to treat models of the form \eqref{eq1}. In the case in which the initial condition decay slowly, \cite{10} states that the solution spreads exponentially faster with a larger index than in the previous case. All these results are in  great contrast with the case $\alpha=1$. They indeed notice that diffusion only plays a role for small times, the large time dynamics being given by a simple transport equation. The scheme of their proof will  be reproduced here, but some steps - and this is why it makes system \eqref{eq1} worth studying - become less easy. The  small time study will require the manipulation of some Polya integrals, and the transport equation will also become more complex.

The work on the single equation \eqref{1} can be extended to  reaction-diffusion systems. The first definitions of spreading speeds
for cooperative  systems in population ecology and epidemic theory are due to Lui in \cite{16}. In a series of papers, Lewis, Li and Weinberger \cite{14},\cite{21}, \cite{19} studied spreading speeds and travelling waves for a particular class of cooperative reaction-diffusion systems,
with standard diffusion. Results on single equations in the singular perturbation framework proved by Evans and Souganidis  in \cite{ES} have also been extended by Barles, Evans and Souganidis in \cite{BES}. The viscosity solutions framework is studied in \cite{BS}, with a precise study of the Harnack inequality. In these papers, the system under study is of the following form

\begin{equation*}
\partial_{t}u_{i} -d_i\Delta u_{i} =f_{i}(u ),
\end{equation*}
where, for   $m\in\N^*$, $u=(u_i)_{i=1}^m$ is the unknown.

For all $i\in\llbracket 1,m\rrbracket:=\{1,...,m\}$, the constants $d_i$ are assumed to be positive as well as the bounded, smooth and Lipschitz initial conditions, defined from $\R^d$ to $\R_+$.
The essential assumptions concern the reaction term $F=(f_i)_{i=1}^m$. This term is assumed to be smooth, to have only two zeroes $0$ and $a\in\R^m$ in $[0,a]$, and for all $i\in\llbracket 1,m\rrbracket$, each $f_i$ is nondecreasing in all its components, with the possible exception of the ith one. The last assumption means that the system   is cooperative. Under additional hypotheses, which imply that the point 0 is unstable, the limiting behaviour of the solution $u=(u_i)_{i=1}^m$ is understood.

Here, we focus on similar systems, keeping the same assumptions on $f$, but considering that at least one diffusive term is given by a fractional Laplacian.  More precisely, we focus on the large time behaviour of the solution $u=(u_{i})_{i=1}^{m}$, for $m\in\mathbb{N}^{*}$,
to the fractional reaction-diffusion system:
\begin{equation}\label{eq1}
\left\{ \begin{array}{rcll}
\partial_{t}u_{i} +(-\triangle)^{\alpha_{i}}u_{i}&=&f_{i}(u ), & t>0, x\in\mathbb{R}^{d},\\
u_{i}(0,x)&=&u_{0i}(x), & x\in\mathbb{R}^{d},
\end{array} \right.
\end{equation}
where
$$\alpha_{i}\in(0,1] \quad \text{ and } \quad \alpha := \displaystyle\min_{\llbracket 1,m \rrbracket} \alpha_i <1.$$

 Note that when $\alpha_{i}=1$, then $(-\triangle)^{\alpha_{i}}=-\triangle$.
 As general assumptions, we impose, for all $i\in\llbracket 1,m \rrbracket $, the initial condition $u_{0i}$ to be nonnegative, non identically equal to 0, continuous  and to satisfy
\begin{equation}\label{eq01} u_{0i}( x)= {\rm{O}}(\abs x^{-(d+2\alpha_{i})}) \quad   \mbox{ as } \  \abs x\rightarrow +\infty. \end{equation}

We also assume that for all $i \in \llbracket 1,m \rrbracket$, the function $f_{i}$ satisfies $f_{i}(0)=0$ and that  system \eqref{eq1} is cooperative, which means:
\begin{equation}\label{f1}f_{i}\in C^{1}(\mathbb{R}^{m}) \ \ \mbox{and} \ \ \partial_j f_{i}>0, \ \text{ on } \R^m, \ \ \mbox{ for all } j \in \llbracket 1,m \rrbracket, \  j\neq i.\end{equation}

We will make additional assumptions on the reaction term $F=(f_i)_{i=1}^m$ that are not general but enable us to understand the long time behaviour of  a class of monotone systems.

The aim of this paper is to understand the time asymptotic location of the level sets of solutions to \eqref{eq1}.
We show that the speed of propagation is exponential in time, with a precise exponent depending on the smallest index $\alpha:=\displaystyle\min_{i\in\llbracket 1,m \rrbracket}\alpha_{i}$
and on the principal eigenvalue of the matrix $DF(0)$ where $F=(f_{i})_{i=1}^{m}$. Also we note that this speed does not depend on the space direction.

For what follows and without loss of generality, we suppose that $\alpha_{i+1}\leq\alpha_{i}$ for
all $i\in\llbracket 1,m-1 \rrbracket$ so that $\alpha=\alpha_{m}<1$. Before stating the main results, we need some additional hypotheses on the nonlinearities $f_{i}$, for all $i\in\llbracket 1,m\rrbracket$.

\begin{itemize}
\item[(H1)] The principal eigenvalue $\lambda_{1}$ of the matrix $DF(0)$ is  positive,
\item[(H2)] There exists $\Lambda>1$ such that, for all $s=(s_i)_{i=1}^m\in \R^m_+$ satisfying $\abs s\geq \Lambda$, we have $f_i(s)\leq 0$,
\item[(H3)] For all $s=(s_i)_{i=1}^m\in\R_+^m$ satisfying  $\abs s\leq \Lambda$, $Df_{i}(0)s-f_{i}(s)\geq c_{\delta_1} {s_{i}}^{1+\delta_{1}}$,
\item[(H4)] For all $s=(s_i)_{i=1}^m\in \R_+^m$ satisfying  $\abs s\leq \Lambda$, $Df_{i}(0)s-f_{i}(s)\leq c_{\delta_2}\abs{s}^{1+\delta_{2}},$
\item[(H5)]  $F=(f_i)_{i=1}^m$ is globally Lipschitz on $\R^m$,
\end{itemize}
where   the constants $c_{\delta_1}$ and $c_{\delta_2}$ are positive and independent of $i\in\llbracket1,m \rrbracket$, and for all $j\in\{1,2\}$
$$\delta_{j}\geq\frac{2}{d+2\alpha}.$$

 This lower bound on $\delta_1$ and $\delta_2$ is a technical assumption to make the supersolution and subsolution to \eqref{eq1}, we construct, to be regular enough. Note that one may easily produce examples of functions $F$ satisfying  (H1) to (H5).

We are now in a position to state our main theorem, which show that the solution to \eqref{eq1} move exponentially fast in time.
\begin{teo}\label{teo1} Let $d\geq1$ and assume that $F$ satisfies (\ref{f1}) and (H1) to (H5). Let $u$ be the solution to \eqref{eq1} with a non negative, non identically equal to 0 and continuous initial condition $u_{0}$ satisfying (\ref{eq01}). Then there exists $\tau>0$ large enough such that for all $i\in\llbracket 1,m \rrbracket$, the following two facts are satisfied:
\begin{itemize}
\item[\textbf{a)}] For every $\mu_{i}>0$, there exists a constant $c>0$ such that,
\[u_{i}(t,x)<\mu_{i}, \ \ \ for \ all \ t\geq\tau \ and \  |x|>ce^{\frac{\lambda_{1}}{d+2\alpha}t}.\]
\item[\textbf{b)}] There exist constants $\varepsilon_{i}>0$ and $C>0$ such that,
\[u_{i}(t,x)>\varepsilon_{i}, \ \ \ for \ all \ t\geq\tau \ and \ |x|<Ce^{\frac{\lambda_{1}}{d+2\alpha}t}.\]
\end{itemize}
\end{teo}

The plan to set  Theorem \ref{teo1} is organized as follows. First, in the short section \ref{mild_solutions}, we state a local existence result of solutions for cooperative systems involving fractional diffusion and  we state a comparison principle for this type of solutions which, although standard, is crucial for the sequel. In Section 3 we deal with finite time and large $x$ decay estimates. The end of this paper, Section 4 is devoted to the proof of Theorem \ref{teo1}.
\setcounter{equation}{0}
\section{Local existence and comparison principle}\label{mild_solutions}
Recall that the operator $A=-{\mathrm{diag}}((-\Delta)^{\alpha_1},\ldots,(-\Delta)^{\alpha_m})$ is sectorial (see \cite{Hy}) in $(L^2(\RR^d))^m$, with
domain $D(A)=H^{2\alpha_1}(\RR^d)\times\ldots\times H^{2\alpha_m}(\RR^d)$. If now $u_0$ satisfies the assumptions of Theorem \ref{teo1}, it is in
$(L^2(\RR^d))^m$, so that the Cauchy Problem \eqref{eq1} has a unique maximal solution, defined on an interval of the form
$[0,t_{max})$; moreover the $L^2$-norm of $u$ blows up as $t\to t_{max}$ if $t_{max}<+\infty$. Finally, we have $u\in C((0,t_{max}),D(A))\cap C([0,t_{max}),(L^2(\RR^d))^m)$ and $\frac{du}{dt}\in C((0,t_{max}),(L^2(\RR^d))^m)$. A standard iteration argument and Sobolev embeddings then yield \[u\in C^p((0,t_{max}),(H^q(\R^d))^m)\] for every integer $p$ and $q$.
\begin{teo}\label{cpm}
Consider $T>0$, and let $u=(u_{i})^{m}_{i=1}$ and $v=(v_{i})^{m}_{i=1}$ such that:   $u\in C((0,T],D(A))\cap C([0,T],(L^2(\mathbb{R}^{d}))^{m})\cap C^{1}((0,T),(L^2(\mathbb{R}^{d}))^{m})
$; and $v\in C([0,T]\times\RR^d)\cap C^1((0,T)\times\RR^d)$. Assume that,
for all $i\in \llbracket 1,m \rrbracket$, we have
\[\partial_{t}u_{i}+(-\triangle)^{\alpha_{i}}u_{i}\leq f_{i}(u), \quad \partial_{t}v_{i}+(-\triangle)^{\alpha_{i}}v_{i}\geq f_{i}(v),\]
where  $f_{i}$  satisfies (\ref{f1}).  If for all $i\in \llbracket 1,m \rrbracket$ and $x\in\mathbb{R}^{d}$, $u_{i}(0,x)\leq v_{i}(0,x)$
we have \[u(t,x)\leq v(t,x) \quad \mbox{ for all } \quad  (t,x)\in[0,T]\times\mathbb{R}^{d}.\]
\end{teo}
{\bf Proof.} Let us define for all $i\in \llbracket 1,m \rrbracket$, $w_{i}=u_{i}-v_{i}$. Then $w_{i}$ satisfies $w_{i}(0,x)\leq0$ and
\begin{eqnarray}
\partial_{t}w_{i}+(-\triangle)^{\alpha_{i}}w_{i}&\leq&f_{i}(u)-f_{i}(v)=\int_{0}^{1}\nabla f_{i}(\sigma u+(1-\sigma)v)d\sigma.(u-v)\nonumber\\
&=&\int_{0}^{1}\nabla f_{i}(\zeta_{\sigma})d\sigma.w,\label{eq9}
\end{eqnarray}
where $\zeta_{\sigma}=\sigma u+(1-\sigma)v$. Notice now that $w_i^+\in C((0,T),H^{2\alpha_i}(\R^d))\cup W^{1,\infty}((0,T),L^2(\mathbb{R}^{d}))$. So,
taking the scalar product of  (\ref{eq9}) with the vector function $(w_{i}^+)_{i=1}^{ m}$ and integrating over $\mathbb{R}^{d}$, we have
\begin{eqnarray}\int_{\mathbb{R}^{d}}w_{i}^{+}\partial_{t}w_{i}dx +\int_{\mathbb{R}^{d}}w_{i}^{+}(-\triangle)^{\alpha_{i}}w_{i}dx\leq\int_{\mathbb{R}^{d}}w_{i}^{+}\int_{0}^{1}\nabla f_{i}(\zeta_{\sigma})d\sigma.w \ dx.\label{c3.8}
\end{eqnarray}
Recall that $\displaystyle\int_{\mathbb{R}^{d}}w_{i}^{+}(-\triangle)^{\alpha_{i}}w_{i}dx\geq 0$. So we have, since $\partial_jf_i(\zeta_\sigma)\geq 0$:
\begin{eqnarray*}
\frac{1}{2}\frac{d}{dt}\left[\int_{\mathbb{R}^{d}}(w_{i}^{+})^{2}dx\right]&\leq& \int_{\mathbb{R}^{d}}\int_{0}^{1}\partial_{i}f_{i}(\zeta_{\sigma})d\sigma(w_{i}^{+})^{2}dx+\sum_{j=1, j\neq i}^{m}\int_{\mathbb{R}^{d}}\int_{0}^{1}\partial_{j}f_{i}(\zeta_{\sigma})d\sigma w_{i}^{+}w_{j}^{+}dx\\
&\leq& C\sum_{j=1}^{m}\int_{\mathbb{R}^{d}}(w_{j}^{+})^{2}dx,
\end{eqnarray*}
where $C$ is a constant that depends on $m$.
Doing this procedure for each $i\in\llbracket 1,m \rrbracket$ and adding, we get for $t \in [0,T]$
\[\frac{d}{dt}\left[\sum_{i=1}^{m}\int_{\mathbb{R}^{d}}(w_{i}^{+})^{2}dx\right]\leq C\sum_{i=1}^{m}\int_{\mathbb{R}^{d}}(w_{i}^{+})^{2}dx.\]
So, by  Gronwall's inequality, we have $w_{i}\leq0$ in $[0,T]\times\RR^d$.
\hfill$\Box$
\setcounter{equation}{0}
\section{Finite time bounds and global existence}
From hypothesis (H2), we deduce that the  positive vector $M=\Lambda\mathsf{1}$, where $\mathsf{1}$ is the vector of size $m$ with all entries equal to $1$, is a supersolution to (\ref{eq1}), if the initial condition $u_{0}=(u_{0i})_{i=1}^m$ is smaller than $ M$. So, from Theorem \ref{cpm}, we have $0\leq u(t,x)\leq M$. To prove global existence, it remains to prove a locally finite $L^2$ bound; this is done in the next subsection.

\subsection{Upper bound}

Now, we are in position to establish an algebraic  upper bound for the solutions of (\ref{eq1}). From (H5), we know that,  for $i\in\llbracket 1,m \rrbracket$ and $j\in\llbracket 1,m \rrbracket$
\[\left|\partial_j f_{i}(s)\right|\leq Lip(f_{i}), \quad \mbox{ for all } s\in\mathbb{R}^{m},\]
where $Lip(f_i)$ is the Lipschitz constant of $f_i$. Taking $l=\max_{i\in\llbracket 1,m \rrbracket}Lip(f_{i})$, we have for all $s=(s_i)_{i=1}^m\geq0$
\begin{equation}\label{eq20}f_{i}(s)=\int_{0}^{1}Df_{i}(\sigma s)d\sigma\cdot s\leq\left|\sum_{j=1}^{m}s_{j}\int_{0}^{1}\frac{\partial f_{i}}{\partial s_{j}}(\sigma s)d\sigma\right|\leq l\sum_{j=1}^{m}s_{j}.\end{equation}

Let us consider $v=(v_{i})^{m}_{i=1}$ the solution of the following system
\begin{equation}\label{eq21} \left\{ \begin{array}{rcll}
\partial_{t}v+Lv&=&Bv,&t>0, x \in \R^m\\
v(0,\cdot)&=&u_{0}, & \R^m,
\end{array} \right. \end{equation}
where $L=\mathrm{diag}((-\triangle)^{\alpha_{1}},...,(-\triangle)^{\alpha_{m}})$, $B=(b_{ij})_{i,j=1}^{m}$ is a matrix with $b_{ij}=l$ for all $i,j\in\llbracket 1,m\rrbracket$. By (\ref{eq20}) and Theorem \ref{cpm}, we conclude that $u\leq v$ in $[0,+\infty)\times\R^d$. A finite time upper bound for $u$ is given by the following lemma.
\begin{lema}\label{lem3}
Let $d\geq1$ and let $u=(u_{i})^{m}_{i=1}$ be the solution of system (\ref{eq1}), with a non negative, non identically equal to 0 and continuous initial condition $u_{0}$ satisfying (\ref{eq01}), and reaction term $F=(f_i)_{i=1}^m$ satisfying \eqref{f1} and (H1) to (H5). Then, for all $i\in\llbracket 1,m \rrbracket$, there exists a locally bounded functions $C_1:(0,+\infty)\rightarrow\mathbb{R}_{+}$ such that for all $t>0$ and $\abs x$ large enough, we have
\[u_{i}(t,x)\leq\frac{C_{1}(t)}{|x|^{d+2\alpha}}.\]
\end{lema}

\noindent
Taking Fourier transforms in each term of system (\ref{eq21}), we have
\begin{equation*}\left\{ \begin{array}{rcll}
\partial_{t}\mathfrak{F}(v)&=&(A(|\xi|)+B)\mathfrak{F}(v),& \xi \in \R^d, t>0\\
\mathfrak{F}(v)(0,\cdot)&=&\mathfrak{F}(u_{0}), & \R^d,
\end{array} \right. \end{equation*}
where $A(|\xi|)=\mathrm{diag}(-|\xi|^{2\alpha_{1}},...,-|\xi|^{2\alpha_{m}})$. Thus, we have that \[\mathfrak{F}(v)(t,\xi)=e^{(A(|\cdot|)+B)t} \ \mathfrak{F}(u_{0})(\xi)\]and then, for all $x\in\R^d$ and $t\geq 0$ :
\begin{equation}\label{sursol}
u(t,x)\leq v(t,x)=\mathfrak{F}^{-1}(e^{(A(|\cdot|)+B)t})\ast u_{0}(x).
\end{equation}
In what follows, we prove that for each time $t>0$, the solution of (\ref{eq1}) decays as $|x|^{-d-2\alpha}$ for large values of $|x|$. Due to the decay of $u_0$ at infinity, we only need to prove that the entries of $\mathfrak{F}^{-1}(e^{(A(|\cdot|)+B)t})$  have the desired decay.
The following lemma is needed to prove that we can rotate the  integration line of a small angle $\eps>0$ in the expression of $\mathfrak{F}^{-1}(e^{(A(|\cdot|)+B)t})$.
\begin{lema}\label{lem0}
For all $z \in \{ z \in \C \ | \ 0 \leq \arg(z) < \frac{\pi}{4\a_1}\}$ and $t\geq 0$, we have
\begin{equation}\label{eq26}\norm{e^{(A(z)+B)t}}\leq e^{(\norm{B}-\abs z^{2\alpha_{1}}\cos(2\alpha_{1}\arg(z)))t}+ e^{(\norm{B}-\abs z^{2\alpha}\cos(2\alpha_{1}\arg(z)))t},
\end{equation}
and if
\begin{equation}\label{It}
I_t(z):=\displaystyle\int_{0}^{t}e^{(t-s)(A(z)+B)}[e^{sB},A(z)]e^{sA(z)}ds,
\end{equation}
then there exists $C_2 : (0,\infty)\rightarrow\mathbb{R}_{+}$  a locally bounded function such that
\begin{equation}\label{estimIt}
\norm{I_t(z)}\leq C_2(t)(\abs z^{2\a}e^{-\abs z^{2\a}\cos(2\a_1\arg(z))t}+\abs z^{2\a_1}e^{-\abs z^{2\a_1}\cos(2\a_1\arg(z))t}).
\end{equation}
\end{lema}
{\bf Proof.} Let $z$ be in $ \{ z \in \C \ | \ 0 \leq \arg(z) < \frac{\pi}{4\a_1}\}$. There exist $j \in \llbracket 1,m\rrbracket$, and $k \in \llbracket 1,m\rrbracket$  such that $\norm{e^{(A(z)+B)t}}= (e^{(A(z)+B)t})_{jk}$. Consider the system
\begin{equation*}\label{eq25}\left\{ \begin{array}{rcll}
\partial_{t}w&=&(A(z)+B)w, & z \in \C, t>0,\\
w(0,z)&=&e_{j}& z \in \C,
\end{array} \right. \end{equation*}
where $e_{j}$ is the $j$th vector of the canonical basis of $\mathbb{R}^{m}$. Thus, we have \[w(t,z)=e^{(A(z)+B)t}.e_{j} \quad \mbox{ and }  \quad \norm{w}= \norm{e^{(A(z)+B)t}}.\]
 Multiply  (\ref{eq25}) by the conjugate transpose $\overline{w}$ and take the real part to get
\[\frac{1}{2}\partial_{t}\norm w^{2}+\sum_{l=1}^{m}\cos(2\alpha_{l}\arg(z))\abs z^{2\alpha_{l}}|w_{l}|^{2}=Re(Bw.\overline{w})\leq \norm B \norm w^{2}.\]
The choice of $\arg(z)$ and Gronwall's Lemma end the proof.

To prove \eqref{estimIt}, it is sufficient to notice that, for $s \in [0,t]$, we have
$$\norm{e^{sA(\abs ze^{i\arg(z)})}}\leq e^{-\abs z^{2\a}\cos(2\a_1\arg(z))s}+e^{-\abs z^{2\a_1}\cos(2\a_1\arg(z))s},
$$
$$
\norm{[e^{sB},A(\abs ze^{i\arg(z)})]}\leq C(t)(\abs z^{2\a}+\abs z^{2\a_1}),
$$
where $C : (0,+\infty)\rightarrow\mathbb{R}_{+}$ is  a locally bounded function, and
due to \eqref{eq26}, we also have
\begin{eqnarray*}
\norm{e^{(A(\abs ze^{i\arg(z)})+B)(t-s)}}  &\leq& e^{(\norm{B}-\abs z^{2\alpha_{1}}\cos(2\alpha_{1}\arg(z)))(t-s)}+e^{(\norm{B}-\abs z^{2\alpha}\cos(2\alpha_{1}\arg(z)))(t-s)}.
\end{eqnarray*} \hfill$\Box$

\noindent
{\bf Proof for $d=1$.} 
In this proof, we denote by  $C:(0,+\infty)\rightarrow\mathbb{R}_{+}$  a locally bounded function. From \eqref{sursol}, we only have to find an upper bound to $\mathfrak{F}^{-1}(e^{(A(|\cdot|)+B)t})$. First, we consider for $t\geq 0$ and $z \in \C$, $w(t,z):=e^{tB}e^{tA(z)}$. Thus, $w$ satisfies the Cauchy problem
\begin{equation*} \left\{ \begin{array}{rcll}
\partial_tw &=&(A(z)+B)w+[e^{tB},A(z)]e^{tA(z)}, & t>0, z \in \C\\
w(0,z)&=&Id,& z \in  \C,
\end{array}\right.
\end{equation*}
where $[e^{tB},A(z)]=e^{tB}A(z)-A(z)e^{tB}$. By Duhamel's formula, we get for all $z \in \C$ and $t\geq 0$:
\begin{eqnarray}\label{exp}e^{t(A(z)+B)}=e^{tB}e^{tA(z)}-\int_{0}^{t}e^{(t-s)(A(z)+B)}[e^{sB},A(z)]e^{sA(z)}ds.\end{eqnarray}

\noindent
Thus, for all $t>0$ and all $x \in \R$, we have
\begin{eqnarray}\label{Fourierinv}
\mathfrak{F}^{-1}(e^{(A(|\cdot|)+B)t})(x)&=&\int_{\mathbb{R}}e^{ix\xi}e^{(A(|\xi|)+B)t}d\xi\\
&=&\int_{\mathbb{R}}e^{ix\xi}e^{tB}e^{tA(\abs \xi)}d\xi-\int_{\mathbb{R}}e^{ix\xi}I_t(\abs \xi)d\xi\notag\\
&=&e^{tB}\ \mathrm{diag}(p_{\a_1}(t,x), ...,p_{\a_m}(t,x)) -\int_{\mathbb{R}}e^{ix\xi}I_t(\abs \xi)d\xi\notag,
\end{eqnarray}
where for $i \in \llbracket 1,m \rrbracket$, $p_{\a_i}$ is the heat kernel of the operator $(-\Delta)^{\a_i}$ in $\R$, \cite{6}. Since $\a=\displaystyle \min_{i\in \llbracket 1,m \rrbracket} \a_i\in(0,1)$,  for large values of $\abs x$, we clearly have
\begin{equation}\label{firstterm}
\norm{e^{tB}\ \mathrm{diag}(p_{\a_1}(t,x), ...,p_{\a_m}(t,x))}\leq \frac{C(t)}{\abs x^{1+2\a}}.
\end{equation}

\noindent
It remains to bound from above the following quantity:
$$\displaystyle\int_{\mathbb{R}}e^{ix\xi}I_t(\abs \xi)d\xi=2\int_{0}^{\infty}\cos(xr)I_t(r)dr=2\re\left( \int_{0}^{\infty}e^{ixr}I_t(r)dr\right).$$

We use the following two facts. First,  for all $t \geq 0$, the function $z \mapsto e^{ixz}I_t(z)$ is holomorphic on $\C \setminus \{0\}$. Second, for $\delta >0$ (respectively $R>0$), on the arc $\{\pm \delta e^{i\theta}, \theta \in [0,\varepsilon]\}$ (respectively $\{ \pm R  e^{i\theta}, \theta \in [0,\varepsilon]\}$), the entries of ${ I_t}$ tends to $0$ as $\delta$ tends to $0$ (respectively $R$ tends to $+\infty$, due to Lemma \ref{lem0}). Consequently, we can rotate the integration line of a small angle $\varepsilon \in (0, \frac{\pi}{4\a_1})$ and the quantity we have to bound from above becomes $\int_{0}^{\infty}e^{ixre^{i\varepsilon}}I_t(re^{i\varepsilon})dr$, with
$$
I_t(re^{i\varepsilon})=\int_{0}^{t}e^{(t-s)(A(re^{i\varepsilon})+B)}[e^{sB},A(re^{i\varepsilon}))]e^{sA(re^{i\varepsilon}))}ds.
$$

\noindent
From Lemma \ref{lem0}, taking $$\eta_{t}=\norm{\int_{0}^{\infty}e^{ixre^{i\varepsilon}}I_t(re^{i\varepsilon})dr}$$  we get,  for  large values of $\abs x$
\begin{eqnarray}\label{secondterm}
\eta_{t}&\leq &C(t) \int_{0}^{\infty}e^{-xr\sin(\varepsilon)}(r^{2\a}e^{-r^{2\a}\cos(2\a_1\varepsilon)t}+r^{2\a_1}e^{-r^{2\a_1}\cos(2\a_1\varepsilon)t})dr\notag\\
&\leq&\frac{C(t)}{\abs x^{1+2\a}} \int_{0}^{\infty}e^{-\tilde r\sin(\varepsilon)}(\tilde r^{2\a} e^{-\frac{\tilde r^{2\a}}{\abs x^{2\a}}\cos(2\a_1\varepsilon)t}+\tilde r^{2\a_1}e^{-\frac{\tilde r^{2\a}}{\abs x^{2\a}}\cos(2\a_1\varepsilon)t})d\tilde r\notag\\
&\leq &\frac{C(t)}{\abs x^{1+2\a}} .
\end{eqnarray}
With \eqref{Fourierinv}, \eqref{firstterm} and \eqref{secondterm}, we conclude that for large values of $\abs x$ and for all $t\geq 0$
$$
\norm{\mathfrak{F}^{-1}(e^{(A(|\cdot|)+B)t})(x)}\leq \frac{C(t)}{\abs x^{1+2\a}},
$$
which concludes the proof.\hfill$\Box$

Now, we state the proof of Lemma \ref{lem3} in the higher space dimension case, i.e. when $d>1$.

\noindent
{\bf Proof.} As previously, from \eqref{sursol}, we only need to bound from above the function $\mathfrak{F}^{-1}(e^{(A(|\cdot|)+B)t})$.
Let $t>0$ and $|x|>1$, using the spherical coordinates system in dimension $d>1$, the definition of Bessel Function of first kind (see \cite{bessel1} and \cite{Erd}), we have
\begin{eqnarray*}
\mathfrak{F}^{-1}(e^{(A(|\cdot|)+B)t})(x)&=&C_{d}\int_{0}^{\infty}\int_{-1}^{1}e^{(A(r)+B)t}\cos(|x|rs)r^{d-1}(1-s^{2})^{\frac{d-3}{2}}dsdr\\
&=&\frac{C_{d}}{|x|^{\frac{d}{2}-1}}\int_{0}^{\infty}e^{(A(r)+B)t}J_{\frac{d}{2}-1}(|x|r)r^{\frac{d}{2}}dr,
\end{eqnarray*}
where $C_d$ is a positive constant depending on $d$.

The matrix $e^{(A(r)+B)t}$ is split into two pieces as done in \eqref{exp}, to get
\begin{eqnarray*}
\mathfrak{F}^{-1}(e^{(A(|\cdot|)+B)t})(x)&=&e^{tB}\ \mathrm{diag}(p_{\a_1}(t,x), ...,p_{\a_m}(t,x))-\frac{C_{d}}{|x|^{\frac{d}{2}-1}}\int_{0}^{\infty}I_t(r)J_{\frac{d}{2}-1}(|x|r)r^{\frac{d}{2}}dr,
\end{eqnarray*}
where $I_t$ has been defined in \eqref{It}.
From  \cite{6}, the first piece of the right hand side has the correct algebraic decay. It remains to bound from above the second piece. In fact, using the Whittaker function (defined in \cite{Erd} for example), we have for all $x\in \R^d$ and all $t>0$:
\begin{eqnarray*}
\frac{C_{d}}{|x|^{\frac{d}{2}-1}}\int_{0}^{\infty}I_t(r)J_{\frac{d}{2}-1}(|x|r)r^{\frac{d}{2}}dr&=&\frac{C_{d}}{|x|^{\frac{d-1}{2}}\sqrt{2\pi}}\re\left(\int_{0}^{\infty}I_t(r)e^{\frac{d-1}{4}i \pi}W_{0,\frac{d}{2}-1}(2i\abs x r)r^{\frac{d-1}{2}}dr\right)\\
&=&\frac{C_{d}}{|x|^{d}\sqrt{2\pi}}\re\left(\int_{0}^{\infty}I_t(\tilde r \abs x^{-1})e^{\frac{d-1}{4}i \pi}W_{0,\frac{d}{2}-1}(2i\tilde r)\tilde r^{\frac{d-1}{2}}d\tilde r\right)
\end{eqnarray*}

As done in the one dimension case, since the Whittaker function is bounded, we can rotate the integration line of a small angle $\varepsilon \in (0, \frac{\pi}{4\a_1})$. Thus, using \eqref{estimIt}, we have the result if we prove that the following integral
$$
\int_{0}^{\infty}\abs{W_{0,\frac{d}{2}-1}(2i\tilde r e^{i\varepsilon})}\tilde r^{\frac{d-1}{2}} (\tilde r^{2\a}+\tilde r^{2\a_1})d\tilde r
$$
is convergent. From \cite{bessel1}, $W_{0,\frac{d}{2}-1}$ has the following asymptotic expressions, thus $W_{0,\frac{d}{2}-1}(z) \underset{\abs z \rightarrow +\infty}{\sim} e^{-\frac{z}{2}}$ and
$$
W_{0,\frac{d}{2}-1}(z) \underset{\abs z \rightarrow 0}{\sim} \begin{cases}
-\Gamma(\frac{d-1}{2})^{-1} \left(\ln(z)+\frac{\Gamma'(\frac{d-1}{2})}{\Gamma(\frac{d-1}{2})}\right)z^{\frac{d-1}{2}}, & \mbox{ if $d=2$} \\
\frac{\Gamma(d-2)}{\Gamma(\frac{d-1}{2})}\  z^{\frac{3-d}{2}},& \mbox{ if $d\geq 3$}.
\end{cases}
$$
\hfill$\Box$
\setcounter{equation}{0}
\subsection{Lower bound}
The following result is   important  and needed to prove   Theorem \ref{teo1}. It sets an algebraically lower bound for the
solutions of the cooperative system (\ref{eq1}). This result is valid for any dimension $d\in \N^*$. Moreover, since for all $i\in\llbracket 1,m\rrbracket$, $f_{i}(0)=0$, we have for all $s=(s_i)_{i=1}^m\in\R^m$ with $0\leq s \leq M$
\[f_{i}(s)=\int_{0}^{1}Df_{i}(\sigma s)d\sigma\cdot s=\sum_{j=1}^{m}s_{j}\int_{0}^{1}\frac{\partial f_{i}}{\partial s_{j}}(\zeta_{\sigma})d\sigma\]
where $\zeta_{\sigma}=\sigma s\in[0,M]$ and $\frac{\partial f_{i}}{\partial s_{j}}:[0,M]\rightarrow\mathbb{R}$ is continuous
for all $i,j\in\llbracket 1,m \rrbracket$, since the system is cooperative, there exist constants $\gamma_{ij}>0$ such that for all $i\in\llbracket 1,m \rrbracket$ and $j\in\llbracket 1,m \rrbracket$:
\begin{equation}\label{gammamm}\left| \partial_if_{i} (\zeta_{\sigma})\right|\leq\gamma_{ii} \ \ \ \mbox{and} \ \ \ \gamma_{ij}\leq \partial_j f_{i} (\zeta_{\sigma}).
\end{equation}

\begin{lema}\label{lem2}
Let $u=(u_{i})^{m}_{i=1}$ be the solution of the system (\ref{eq1}), with non negative, non identically equal to 0 and continuous initial condition $u_{0}$ satisfying (\ref{eq01}) and with reaction term $F=(f_i)_{i=1}^m$ satisfying   (\ref{f1}), (H1), (H2) and (H5).
Then, for all $i\in \llbracket 1,m\rrbracket$,  $x\in\mathbb{R}^{d}$ and $t\geq{1}$, we have:
\begin{equation}\label{lowerbound} u_{i}(t ,x)\geq\frac{\underline c \ t \ e^{-\gamma_{mm}t}}{t^{\frac{d}{2\alpha}+1}+|x|^{d+2\alpha}},
\end{equation}
 where $\underline c$ is a positive constant and $\gamma_{mm}$ is defined in \eqref{gammamm}.
\end{lema}

\noindent
{\bf Proof.} We split the proof into three steps: first,  we prove the result  for $i=m$, which serves as an initiation of the process. In an intermediate step,
for all $i\in\llbracket 1,m-1 \rrbracket$, $t\geq 1$ and $s \in [0, t-1]$, we find a lower bound of $p_{\a_i}(\cdot,t-s)\star (s^{\frac{d}{2\a}+1}+\abs \cdot ^{d+2\a})^{-1}$, that decays like $\abs x^{-(d+2\a)}$ for large values of $\abs x$. In a third step,  for all $i\in\llbracket 1,m-1 \rrbracket$, $t\geq 1$ and $s \in [0, t-1]$, we prove that $u_i(t,\cdot)$ can be bounded from below by  an expression that only depends on the integral $\displaystyle\int_0^t p_{\a_i}(\cdot,t-s)\star (s^{\frac{d}{2\a}+1}+\abs \cdot ^{d+2\a})^{-1}ds$.

\noindent
\textsf{Step 1.} We have for all $x\in\mathbb{R}^{d}$ and $t>0$:
\[\partial_{t}u_{m}+(-\triangle)^{\alpha_{m}}u_{m}=f_{m}(u)\geq \int_{0}^{1} \partial_m f_{m} (\zeta_{\sigma})d\sigma u_{m}\geq-\gamma_{mm}u_{m},\]
where $\gamma_{mm}$ is defined in \eqref{gammamm}.
By the maximum principle of reaction diffusion equations, we have for all $ t\geq0$
\[u_{m}(t,x)\geq e^{-\gamma_{mm}t}(p_{\a_m}(t,\cdot)\ast u_{0m})(x), \]
Since $u_{0m}(\cdot)\not\equiv0$ is continuous and nonnegative, we can find $\xi\in\mathbb{R}^{d}$ such that $u_{0m}(y)\geq C$ for all $y\in B_{R}(\xi)$ for some $R>0$ and $C>0$. If $|x|>R$, $t\geq1$ and using that $\alpha:=\alpha_{m}<1$, we get
\begin{eqnarray*}
(p_{ \a_m}(t,\cdot)\ast u_{0m})(x)&\geq&C\int_{\mid y-\xi\mid\leq R}\frac{B^{-1}t}{t^{\frac{d}{2\alpha}+1}+|x-y|^{d+2\alpha}}dy\\
&=&C\int_{|z|\leq R}\frac{B^{-1}t}{t^{\frac{d}{2\alpha}+1}+|x-\xi-z|^{d+2\alpha}}dz.
\end{eqnarray*}
We also have $|x-\xi-z|\leq(2+\frac{\xi}{R})|x|$. Thus
\[t^{\frac{d}{2\alpha}+1}+|x-\xi-z|^{d+2\alpha}\leq \left(2+\frac{\xi}{R}\right)^{d+2\alpha}t^{\frac{d}{2\alpha}+1}+\left(2+\frac{\xi}{R}\right)^{d+2\alpha}|x|^{d+2\alpha}.\]
Then
\begin{eqnarray*}
(p_{\a_m}(t,\cdot)\ast u_{0m})(x)&\geq&\frac{CB^{-1}}{(2+\frac{\xi}{R})^{d+2\alpha}}\int_{|z|\leq R}\frac{t}{t^{\frac{d}{2\alpha}+1}+|x|^{d+2\alpha}}dz=\frac{\widetilde{C}t}{t^{\frac{d}{2\alpha}+1}+|x|^{d+2\alpha}},
\end{eqnarray*}
where $\tilde C$ is a positive constant.
If $|x|\leq R$ and $t\geq1$,
\begin{eqnarray*}
(p_{\a_m}(t,\cdot)\ast u_{0m})(x)&\geq&\int_{B_{1}(0)}\frac{B^{-1}t}{t^{\frac{d}{2\alpha}+1}+|x-y|^{d+2\alpha}}u_{0m}(y)dy\\
&\geq&\frac{B^{-1}t}{t^{\frac{d}{2\alpha}+1}+(R+1)^{d+2\alpha}}\int_{B_{1}(0)}u_{0m}( y)dy\\
&\geq&\frac{\overline{C}t}{t^{\frac{d}{2\alpha}+1}} \geq\frac{\overline{C}t}{t^{\frac{d}{2\alpha}+1}+|x|^{d+2\alpha}} ,
\end{eqnarray*}
for some small constant $\overline{C}>0$.
Then, there exist $C_{m}>0$ such that  for all $x\in\mathbb{R}^{d}$ and $t\geq1$
\begin{equation}\label{sousestimu}u_{m}(t,x)\geq\frac{C_{m}te^{-\gamma_{mm}t}}{t^{\frac{d}{2\alpha}+1}+|x|^{d+2\alpha}}.
\end{equation}
\ \\\textsf{Step 2.} By similar computations as done in Step 1, it is possible to find a constant $C>0$ such that for all $x\in\mathbb{R}^{d}$,  $t\geq1$ and $s\in[0,t-1]$:
\begin{itemize}
\item[-] if $\a_i=1$ then
\begin{eqnarray*}
p_{\a_i}(\cdot,t-s)\star (s^{\frac{d}{2\a}+1}+\abs \cdot ^{d+2\a})^{-1}(x)&\geq&\frac{1}{(4\pi (t-s))^{\frac{d}{2}}}\int_{\R^d}\frac{e^{-\frac{|y|^{2}}{4(t-s)}}}{s^{\frac{d}{2\alpha}+1}+|x-y|^{d+2\alpha}}dy \\
&\geq&\frac{1}{(4\pi (t-s))^{\frac{d}{2}}(s^{\frac{d}{2\alpha}+1}+|x |^{d+2\alpha})},
\end{eqnarray*}
\item[-] if $\a_i \in (0,1)$ then
\begin{eqnarray*}
p_{\a_i}(\cdot,t-s)\star (s^{\frac{d}{2\a}+1}+\abs \cdot ^{d+2\a})^{-1}(x)&\geq&\int_{\R^d}\frac{1}{((t-s)^{\frac{d}{2\alpha_i}+1}+\abs y^{d+2\alpha_i})(s^{\frac{d}{2\alpha}+1}+|x-y|^{d+2\alpha})}dy\notag\\
&\geq&\frac{(t-s)^{-\frac{d}{2\a_i}}}{s^{\frac{d}{2\alpha}+1}+\abs x^{d+2\alpha}}.
\end{eqnarray*}
\end{itemize}
\ \\\textsf{Step 3.} For $i\in\llbracket 1,m-1 \rrbracket$, we have for all $x\in\mathbb{R}^{d}$ and $t\geq 0$
\begin{eqnarray*}
\partial_{t}u_{i}+(-\triangle)^{\alpha_{i}}u_{i}\geq \int_{0}^{1} \partial_m f_{i}(\zeta_{\sigma})d\sigma u_{m}+\int_{0}^{1}\partial_i f_{i}(\zeta_{\sigma})d\sigma u_{i}\geq\gamma_{im}u_{m}-\delta_{i}u_{i},
\end{eqnarray*}
where $\zeta_{\sigma}= \sigma u$ and $\delta_{i}\geq\max(\gamma_{ii},\gamma_{m}+1)$. Then,
by the maximum principle of reaction diffusion equations and Duhamel's formula, we have for all $(t,x)\in\mathbb{R}_{+}\times\mathbb{R}^{d}$
\begin{eqnarray*}
u_{i}(t,x)\geq e^{-\delta_{i}t}(p_{\a_i}(t,\cdot)\ast u_{0i})(x)+\gamma_{im}e^{-\delta_{i}t}\int_{0}^{t}\int_{\mathbb{R}^{d}}p_{\a_i}(t-s,y)u_{0m}(s,x-y)e^{\delta_{i}s}dyds.
\end{eqnarray*}
So, taking $t\geq1$, and using \eqref{sousestimu}, we get
\begin{eqnarray*}
u_{i}(t,x)
\geq C_{m}\gamma_{im}e^{-\delta_{i}t}\int_{0}^{t-1}\int_{\mathbb{R}^{d}}p_{\a_i}(t-s,y)\frac{se^{(\delta_{i}-\gamma_{mm})s}}{s^{\frac{d}{2\alpha}+1}+|x-y|^{d+2\alpha}}dyds
\end{eqnarray*}
Using Step 2, we get the following lower bound, for all $x\in\mathbb{R}^{d}$, $t\geq1$, taking $C_i$ smaller if necessary:
\begin{eqnarray*}
u_{i}(t,x)&\geq& C_{i}\frac{e^{-\delta_{i}t}}{t^{\frac{d}{2\alpha}}}\int_{0}^{t-1}\frac{se^{(\delta_{i}-\gamma_{mm})s}-1)s}{s^{\frac{d}{2\alpha}+1}+|x|^{d+2\alpha}}ds\geq\frac{C_{i}te^{-\gamma_{mm}t}}{t^{\frac{d}{2\alpha}+1}+|x|^{d+2\alpha}}.
\end{eqnarray*}
\hfill$\Box$
\setcounter{equation}{0}
\section{Proof of Theorem \ref{teo1}}
Inspired by the formal analysis done in \cite{7}, we construct an explicit supersolution (respectively subsolution) of the form
\begin{equation}\label{eq2} v(t,x)=a\left(1+b(t)|x|^{\delta(d+2\alpha)}\right)^{-\frac{1}{\delta}}\phi_1,\end{equation}
where $b(t)$ is a time continuous function  asymptotically proportional to $e^{-\delta\lambda_{1} t}$,  $\phi_1=(\phi_{1,i})_{i=1}^m\in\mathbb{R}^{m}$ is the normalised (positive)
principal eigenvector of $DF(0)$ associated to the principal eigenvalue $\lambda_{1}$, and $\delta$ is equal to $\delta_1$ (respectively $\delta_2$) defined in (H3) (respectively (H4)).

\begin{lema}\label{lem1} Let $v$ be defined as in (\ref{eq2}). Then,   there exist a constant $D>0$ such that for all $i\in\llbracket 1,m \rrbracket$, $t>0$ and $x \in \R^d$
\[\mid(-\triangle)^{\alpha_{i}}v_{i}(t,x)\mid\leq Db(t)^{\frac{2\alpha_{i}}{\delta(d+2\alpha)}}v_{i}(t,x),\]
where $\alpha_{i}\in(0,1]$.
\end{lema}
{\bf Proof.} The case $\alpha_{i}=1$ is trivial.
For $\alpha_{i}\in(0,1)$ and $\delta\geq\displaystyle\frac{2}{d+2\alpha}$, since $(-\Delta)^{\a_i}$ is $2\alpha_{i}$-homogeneous, we only need to prove
$$
\mid(-\triangle)^{\alpha_{i}} w(x)\mid\leq Dw(x)
$$
where $w(x)=(1+| x|^{\delta(d+2\alpha)})^{-\frac{1}{\delta}}$.

We consider the following decomposition, which is the central part of the proof:
\begin{eqnarray*}
(-\triangle)^{\alpha_{i}} w(x)&=&\int_{|y|>3|x|/2}\frac{w(x)-w(y)}{|x-y|^{d+2\alpha_{i}}}dy+\int_{B_{|x|/2}(x)}\frac{w(x)-w(y)}{|x-y|^{d+2\alpha_{i}}}dy\\
&&\quad\quad+\int_{\{|x|\leq2|y|\leq3|x|\}\setminus B_{|x|/2}(x)}\frac{w(x)-w(y)}{|x-y|^{d+2\alpha_{i}}}dy+\int_{|y|\leq|x|/2}\frac{w(x)-w(y)}{|x-y|^{d+2\alpha_{i}}}dy.
\end{eqnarray*}
Each piece is easily bounded, as in  \cite{4} for instance. \hfill$\Box$

\medskip
In what follows, we will use the results of previous sections to obtain appropriate sub and super solutions to (\ref{eq1}) of the form (\ref{eq2}).
We divide the proof of Theorem \ref{teo1} in two lemmas.

\begin{lema}\label{lem1.teo1} Assume that $F$ satisfies (\ref{f1}), (H1), (H2), (H3) and (H5). Let $u$ be the solution to (\ref{eq1}) with $u_{0}$ satisfying the assumptions of Theorem \ref{teo1}. Then, for every $\mu=(\mu_{i})_{i=1}^{m}>0$, there exists $c>0$ such that, for all $t>\tau$, with $\tau>0$ large enough

\[\left\{x\in\mathbb{R}^{d}\mid \ |x|>ce^{\frac{\lambda_{1}}{d+2\alpha}t}\right\}\subset\left\{x\in\mathbb{R}^{d}\mid \ u(t,x)<\mu\right\}.\]
\end{lema}
{\bf Proof:} We consider the function $\overline{u}$ given by (\ref{eq2}) with $\delta=\delta_{1}$ as in (H3). The idea is to adjust $a>0$ and $b(t)$ so that the function
$\overline{u}$ serves as supersolution of \eqref{eq1}.

 In the sequel,  $a$ is any positive constant satisfying \[a\geq\left(\frac{D+\lambda_{1}}{c_h} \right)^{\frac{1}{\delta_{1}}}\displaystyle \max_{i\in\llbracket1,m\rrbracket}\left(\frac{1}{\phi_{1,i}}\right),\]
where $c_h$ is defined in (H2).
For any constant $B\in(0,(1+D\lambda_{1}^{-1})^{-\frac{\delta_{1}(d+2\alpha)}{2\alpha}})$, where $D>0$ is given in Lemma \ref{lem1}, we consider the following ordinary differential equation
\begin{equation}\label{equationb}b'(t)+\delta_{1}Db(t)^{\frac{2\alpha}{\delta_{1}(d+2\alpha)}+1}+\delta_{1}\lambda_{1}b(t)=0, \ \ b(0)=(-D\lambda_{1}^{-1}+B^{-\frac{2\alpha}{\delta_{1}(d+2\alpha)}})^{-\frac{\delta_{1}(d+2\alpha)}{2\alpha}}
\end{equation}
whose solution is given by
$$b(t)=(-D\lambda_{1}^{-1}+B^{-\frac{2\alpha}{\delta_{1}(d+2\alpha)}}e^{\frac{2\alpha\lambda_{1}}{d+2\alpha}t})^{-\frac{\delta_{1}(d+2\alpha)}{2\alpha}}.
$$
For all $t \geq 0$, we have $0\leq b(t)\leq b(0)\leq1$. Using Lemma \ref{lem1}, we have for all $i\in\llbracket 1,m \rrbracket$
$$
\begin{array}{rll}
&&\partial_{t}\overline{u}_{i}+(-\triangle)^{\alpha_{i}}\overline{u}_{i}-f_{i}(\overline{u})=\partial_{t}\overline{u}_{i}+(-\triangle)^{\alpha_{i}}\overline{u}_{i}-Df_{i}(0)\overline{u}+\left[Df_{i}(0)\overline{u}-f_{i}(\overline{u})\right]\quad\\
&&\quad\quad\quad\quad\quad\geq\displaystyle\frac{a\phi_{1,i}}{\delta_{1}(1+b(t)|x|^{\delta_{1}(d+2\alpha)})^{\frac{1}{\delta_{1}}+1}}\left\{-b'(t)-\delta_{1}Db(t)^{\frac{2\alpha}{\delta_{1}(d+2\alpha)}+1}-\delta_{1}\lambda_{1}b(t)\right\}|x|^{\delta_{1}(d+2\alpha)}\nonumber\\
&&\quad\quad\quad\quad\quad\quad\quad\quad+\displaystyle\frac{a\phi_{1,i}}{(1+b(t)|x|^{\delta_{1}(d+2\alpha)})^{\frac{1}{\delta_{1}}+1}}\left\{-Db(t)^{\frac{2\alpha}{\delta_{1}(d+2\alpha)}}-\lambda_{1}+c_h\phi_{1,i}^{\delta_{1}}a^{\delta_1}\right\}\geq 0.
\end{array}
$$
Finally, due to  Lemma \ref{lem3},  for a fixed $t_{0}>0$, there exists $t_{1}\geq 0$ such that for all $x\in\mathbb{R}^{d}$ and all $ i\in\llbracket 1,m \rrbracket$, we have $\overline{u}_{i}(t_{1},x)\geq u_{i}(t_{0},x).$
Thus, by Theorem \ref{cpm} we have,  for all $t\geq t_{0}$, all $x\in\mathbb{R}^{d}$ and all $i\in\llbracket 1,m \rrbracket$ :
$\overline{u}_{i}(t+t_{1}-t_{0},x)\geq u_{i}(t,x).$

\noindent
For  any $(\mu_{i})_{i=1}^{m}>0$, we define for $i\in\llbracket1,m\rrbracket$ the constants
$$c_{i}^{d+2\alpha}:=a\phi_{1,i}e^{\lambda_{1}(t_{1}-t_{0})}[\mu_{i}B^{\frac{1}{\delta_{1}}}]^{-1}.
$$
Taking $c=\displaystyle \max_{i\in\llbracket1,m\rrbracket} c_i$, if $|x|>c e^{\frac{\lambda_{1}}{d+2\alpha}t}$, then,
for all $t>\tau:=t_{0}$ and all $i\in\llbracket 1,m \rrbracket$
$$ u_{i}(t,x)\leq\overline{u}_{i}(t+t_{1}-t_{0},x)=\frac{a\phi_{1,i}}{(1+b(t+t_{1}-t_{0})|x|^{\delta_{1}(d+2\alpha)})^{\frac{1}{\delta_{1}}}}<\mu_{i}.
$$
\hfill$\Box$

\begin{lema}\label{lem2.teo1} Let $d\geq1$ and assume that $F$ satisfies \eqref{f1}, (H1), (H2), (H4) and (H5). Let $u$ be the solution to (\ref{eq1}) with a non negative, non identically equal to 0 and continuous initial condition $u_{0}$ satisfying (\ref{eq01}).
Then, for all $i\in\llbracket 1,m \rrbracket$, there exist constants $\varepsilon_{i}>0$ and $C>0$ such that,
\[u_{i}(t,x)>\varepsilon_{i}, \quad  \mbox{for  all } \quad t\geq t_1 \ \mbox{ and  } \ |x|<Ce^{\frac{\lambda_{1}}{d+2\alpha}t},\]
with $t_1>0$ large enough.
\end{lema}

{\bf  Proof:} As in the previous proof, we consider the function $\underline{u}$ given by (\ref{eq2}) with $\delta=\delta_{2}$ defined in (H4).
Since, $\underline{u}_{i}(0,\cdot)\leq u_{0i}$ may not hold for all $i\in\llbracket 1,m \rrbracket$, we look for a time $t_{1}>0$ such that $\underline{u}_{i}(0,\cdot)\leq u_{i}(t_{1},\cdot)$ for all $i\in\llbracket 1,m \rrbracket$. Indeed, let $L$ be a constant greater than $\max\{D,\lambda_{1}\}$, where $D$ is given by Lemma \ref{lem1}.
We choose $t_1\geq \max(1,2D\lambda_1^{-1})$ large enough, so that if we set
\begin{equation}\label{c3.10}
a=\frac{\displaystyle\min_{i\in\llbracket1,m\rrbracket}C_{i}\ e^{-\gamma_{mm}t_{1}}}{2\displaystyle\max_{i\in\llbracket1,m\rrbracket}\phi_{1,i}\ t_{1}^{\frac{d}{2\alpha}}} \quad  \mbox{ and }  \quad B=\left(\frac{2}{t_{1}}\right)^{\frac{(d+2\alpha)}{2\alpha}\delta_{2}},
\end{equation}
then
$$
a\leq \left(\frac{\displaystyle\min_{i\in\llbracket1,m\rrbracket}\phi_{1,i}\ \lambda_{1}}{2c_{\delta_2}}\right)^{\frac{1}{\delta_{2}}} \quad \mbox{ and }  \quad B \leq (D\lambda_{1}^{-1})^{-\frac{(d+2\alpha)}{2\alpha}\delta_{2}},
$$
where $c_{\delta_2}$ is defined in (H4).
Then we set
\begin{eqnarray*}b(t)=(D\lambda_{1}^{-1}+B^{-\frac{2\alpha}{\delta_{2}(d+2\alpha)}}e^{\frac{2\alpha\lambda_{1}}{d+2\alpha}t})^{-\frac{(d+2\alpha)}{2\alpha}\delta_{2}}. \end{eqnarray*}
Using Lemma \ref{lem1} and (H3), similarly to the previous proof, we can state that, $\mbox{for all} \ i\in\llbracket 1,m \rrbracket,$
\begin{eqnarray}
\partial_{t}\underline{u}_{i}+(-\triangle)^{\alpha_{i}}\underline{u}_{i}-f_{i}(\underline{u})\leq0, \quad \quad \mbox{in} \ (0,+\infty)\times \R^d.\nonumber
\end{eqnarray}

\noindent
From Lemma \ref{lem2}, we know that for all $i\in\llbracket 1,m \rrbracket$ and all $ x\in\mathbb{R}^{d}$
\[u_{i}(t_{1},x)\geq\underline c\frac{t_{1}e^{-\gamma_{mm}t_{1}}}{t_{1}^{\frac{d}{2\alpha}+1}+|x|^{d+2\alpha}}.\]
 By (\ref{c3.10}),
we deduce
\begin{eqnarray*}\underline c t_{1}e^{-\gamma_{mm}t_{1}}(1+b(0)|x|^{\delta_{2}(d+2\alpha)})^{\frac{1}{\delta_{2}}}\geq\frac{\underline c}{2}t_{1}e^{-\gamma_{mm}t_{1}}(1+b(0)^{\frac{1}{\delta_{2}}}|x|^{d+2\alpha})\geq a\phi_{i}(t_{1}^{\frac{d}{2\alpha}+1}+|x|^{d+2\alpha}).
\end{eqnarray*}
Therefore, we get, for all $i\in\llbracket 1,m \rrbracket$, $u_{i}(t_{1},\cdot)\geq \underline{u}_{i}(0,\cdot)$ in $\mathbb{R}^{d}$, and by Theorem \ref{cpm}, we have for all $t \geq t_1$
\[u_{i}(t,\cdot)\geq\underline{u}_{i}(t-t_{1},\cdot), \quad\quad \mbox{ in } \R^d\]
Finally we choose
\[ \varepsilon_{i}=\frac{a\phi_{1,i}}{2^{\frac{1}{\delta_{2}}}} \quad \mbox{ and } \quad C^{d+2\alpha}=e^{-\lambda_{1}t_{1}}B^{-\frac{1}{\delta_{2}}}.\]
If $t\geq\tau:=t_{1}$ and $|x|\leq Ce^{\frac{\lambda_{1}}{d+2\alpha}t}$, we have
\[u_{i}(t,x)\geq \underline{u}_{i}(t-t_{1},x)=\frac{a\phi_{1,i}}{(1+b(t-t_{1})|x|^{\delta_{2}(d+2\alpha)})^{\frac{1}{\delta_{2}}}}\geq\frac{a\phi_{1,i}}{2^{\frac{1}{\delta_{2}}}}=\varepsilon_{i}.\]


\begin{thebibliography}{AAA}
\bibitem{17} A.N. Kolmogorov, I.G. Petrovsky and N.S. Piskunov, \emph{\'Etude de l'\'equation de la diffusion avec croissance de la quantit\'e de mati\`ere et son application \`a un probl\`eme biologique}, Bull. Univ. \'Etat Moscou S\'er. Inter. A 1 (1937) 1-26.
\bibitem{1} D.G. Aronson and H.F. Weinberger, \emph{Multidimensional nonlinear diffusions arising in population genetics}, Adv. Math. 30 (1978), 33-76.
\bibitem{6} X. Cabr\'e and J. Roquejoffre. \emph{The influence of fractional diffusion in Fisher-KPP equation}. Preprint, arXiv:1202.6072v1, (2012).
\bibitem{7} X. Cabr\'e, A. C. Coulon, and J. M. Roquejoffre. \emph{Propagation in Fisher-KPP type equations with fractional diffusion in periodic media}. C. R. Math. Acad. Sci. Paris, 350 (2012), no. 19-20, 885-890.
\bibitem{10} P. Felmer and M. Yangari. \emph{Fast Propagation for Fractional KPP Equations with Slowly Decaying Initial Conditions}. SIAM J. Math. Anal., 45(2), 662-678.
\bibitem{16} R. Lui. \emph{Biological growth and spread modeled by systems of recursions. I. Mathematical theory}. Math. Biosci. 93(2), 269-295 (1989).
\bibitem{14} M. Lewis, B. Li and H. Weinberger. \emph{Spreading speed and linear determinacy for two-species competition models}. J. Math. Biol. 45, 219-233 (2002).
\bibitem{21} H.F. Weinberger, M. Lewis and B. Li. \emph{Anomalous spreading speeds of cooperative recursion systems}. J. Math. Biol. 55, 207-222 (2007).
\bibitem{19} H.F. Weinberger, M. Lewis and B. Li. \emph{Analysis of linear determinacy for spread in cooperative models}. J. Math. Biol. 45, 183-218 (2002).
\bibitem{ES} L. C. Evans  and P. E. Souganidis, \emph{A {PDE} approach to geometric optics for certain semilinear  parabolic equations}, Indiana Univ. Math. J. 45(2) (1989), 141--172.
\bibitem{BES} G. Barles and L. C. Evans and P. E. Souganidis, \emph{Wavefront propagation for reaction-diffusion systems of {PDE}}. Duke Math. J. 61 (1990), 835-858.
\bibitem{BS} J. Busca and B. Sirakov, \emph{Harnack type estimates for nonlinear elliptic systems and applications}, Ann. Inst. H. Poincar\'e Anal. Non Lin\'eaire 21 (2004), 543--590.
\bibitem{Hy} D. Henry, \emph{Geometric Theory of Semilinear Parabolic Equations}, Springer-Verlag, New York (1981).
\bibitem{bessel1} M. Abramowitz and I. A. Stegun,  \emph{Handbook of Mathematical Functions With Formulas, Graphs, and Mathematical Tables}. Dover Publications, New York. 1972.
\bibitem{Erd}A. Erd{\'e}lyi, \emph{Higher Transcendental Functions. {V}ol. {I}}, New York-Toronto-London, McGraw-Hill Book Company, Inc.,1953.
\bibitem{4} M. Bonforte and J. Vazquez. \emph{Quantitative Local and Global A Priori Estimates for Fractional Nonlinear Diffusion Equations}. Preprint, ArXiv:1210.2594.

\end{thebibliography}
\end{document}